\documentclass{amsart}

\usepackage{amsmath}
\usepackage{amsthm}
\usepackage{amssymb}
\usepackage[all]{xy}
\usepackage[mathscr]{eucal}
\usepackage[german,french,english]{babel}

\UseComputerModernTips

\newcommand{\nat}{\ensuremath{\mathbb{N}}}
\newcommand{\rat}{\ensuremath{\mathbb{Q}}}
\newcommand{\ganz}{\ensuremath{\mathbb{Z}}}

\newcommand{\C}{\ensuremath{\mathbb{C}}}

\newcommand{\eins}{1 \hspace{-2.3pt} \mathrm{l}}

\newcommand{\Knull}{\ensuremath{\operatorname{K}_0}}

\newcommand{\Vark}{\ensuremath{\mathit{Var}_k}}
\newcommand{\Bl}[2]{\ensuremath{\operatorname{Bl}_{#1}#2}}
\newcommand{\cA}{\ensuremath{\mathcal{A}}}

\newcommand{\cM}{\ensuremath{\mathcal{M}}}
\newcommand{\tate}{\ensuremath{\mathbb{L}}}
\newcommand{\aff}[1]{\ensuremath{\mathbb{A}^{#1}}}
\newcommand{\proj}[1]{\ensuremath{\mathbb{P}^{#1}}}
\newcommand{\Proj}{\ensuremath{\operatorname{\mathbb{P}}}}

\newcommand{\Sign}{\ensuremath{\operatorname{Sign}}}
\newcommand{\Triv}{\ensuremath{\operatorname{Triv}}}
\newcommand{\End}{\ensuremath{\operatorname{End}}}
\newcommand{\Hom}{\ensuremath{\operatorname{Hom}}}
\newcommand{\cHom}{\ensuremath{\operatorname{\mathscr{H}\hspace{-1pt}\mathit{om}}}}
\newcommand{\Ind}{\ensuremath{\operatorname{Ind}}}
\newcommand{\Rest}{\ensuremath{\operatorname{Res}}}
\newcommand{\Sym}{\ensuremath{\operatorname{Sym}}}
\newcommand{\Alt}{\ensuremath{\operatorname{Alt}}}

\newcommand{\Pic}{\ensuremath{\operatorname{Pic}}}
\newcommand{\im}{\ensuremath{\operatorname{im}}}
\newcommand{\Chowk}{\ensuremath{\Knull(\mathit{CM}_k)}}

\newtheorem{thm}{Theorem}[section]

\newtheorem{lem}[thm]{Lemma}
\newtheorem{prop}[thm]{Proposition}

\theoremstyle{definition}

\theoremstyle{remark}
\newtheorem{rem}[thm]{Remark}

\numberwithin{equation}{section}

\begin{document}
\title[Functional equations for motivic zeta functions]{A note on functional equations for 
zeta functions with values in Chow motives}
\author{Franziska Heinloth}
\address{Universit\"at Duisburg---Essen, Standort Essen, 
FB6, Mathematik, 45117 Essen, Germany}
\email{franziska.heinloth@uni-duisburg-essen.de}
\begin{abstract}
We consider zeta functions with values in the Grothendieck ring of
Chow motives. Investigating the $\lambda$--structure of this ring,
we deduce a functional equation for the zeta function of abelian varieties.
Furthermore, we show that the property of
having a rational zeta function satisfying a functional equation is preserved
under products.
\end{abstract}
\maketitle

\section{Introduction}
Let $C$ be a geometrically irreducible
smooth projective curve of genus $g$ over a field $k$.
Kapranov in \cite{Kapranov} considers the \emph{zeta function} 
$$\zeta_\mu(C,T) =\sum_{i=0}^{\infty} \mu(\Sym^i(C))T^i,$$
where $\mu$ is a multiplicative Euler characteristic with compact support (i.e. an 
invariant of $k$-varieties with values in a ring $A$ 
satisfying $\mu(X)=\mu(X-Y)+\mu(Y)$ for $Y\subset X$ closed and 
$\mu(X\times Y)=\mu(X)\times\mu(Y)$), and $\Sym^i(C)$ denotes the $i$-th
symmetric power of $C$. For example, if $k$ is a finite field, the
number of $k$-valued points is such an invariant, and the associated zeta 
function is the Hasse--Weil zeta function.
Kapranov shows that if $A$ is a field and $\tate_\mu=\mu(\aff{1})\ne0$, the zeta
function of $C$ with respect to $\mu$ is rational and satisfies the functional equation
$$\zeta_\mu(C,\frac{1}{\tate_\mu T})=\tate_\mu^{1-g}T^{2-2g}\zeta_\mu(C,T).$$

Kapranov suggests that also zeta functions of higher dimensional smooth
projective varieties should be rational and satisfy a functional equation.

Larsen and Lunts in \cite{LarsenLunts} and \cite{LarsenLunts2} 
for $k=\C$ construct  
a multiplicative Euler characteristic with
compact support $\mu$ 
such that the zeta function with respect to $\mu$
of smooth projective surfaces of nonnegative Kodaira dimension
is not rational. (In their example, 
$\tate_\mu=0$.)

On the other hand, as they point out in
\cite{LarsenLunts2}, if $A$ carries a $\lambda$--structure $\sigma^i$ such
that $A$ (with its opposite structure) is 
\emph{special} (compare Section \ref{lambda}), and if $\mu(\Sym^i X)=\sigma^i(\mu(X))$,
the property of having a
rational zeta function is e.g. preserved under products.

In this note, we consider the value ring $\Chowk$, the Grothendieck ring of 
Chow motives over $k$ with rational coefficients.
It is the free abelian group on isomorphism classes $[M]$ of Chow motives $M$
modulo the relations $[M\oplus N]=[M]+[N]$ and carries a commutative ring
structure induced by the tensor product of Chow motives.
There is also the notion of the $i$-th symmetric power $\Sym^i M$ 
of a Chow motive
$M$, which is defined as the image of the projector
$\frac{1}{i!}\sum_{\sigma\in S_i}\sigma$ on $M^{\otimes i}$. 
The symmetric powers $\Sym^i$ endow $\Chowk$
with the structure of a $\lambda$-ring. The opposite structure
$(\Alt^i)_i$ is induced by the 
projectors
$\frac{1}{i!}\sum_{\sigma\in S_i} (-1)^\sigma \sigma$
 and turns out to
be special (see Section \ref{lambdaknull} for details).

In characteristic zero, Gillet and Soul\'e as a corollary from
\cite{GilletSoule} and Guillen and Navarro Aznar as a corollary from
\cite{GuillenNavarro}
get a multiplicative Euler characteristic with compact support $\mu$ with
values in $\Chowk$, such that 
$\mu(X)=[h(X)]$ for a smooth projective variety $X$. Here $h(X)$ is the Chow 
motive of $X$ . Note that $\mu(\aff{1})$ is the class of the
Tate motive $\tate$.
It follows from a result of Del Ba{\~n}o and Navarro Aznar in 
\cite{BanoNavarro} that
that $\mu(\Sym^i X)=[\Sym^i h(X)]$ for a smooth projective variety $X$. 
Hence the zeta function of $X$ associated to $\mu$
equals 
$$Z_X(T)=\sum_{i=0}^{\infty} [\Sym^ih(X)]T^i.$$

This zeta function with values in $\Chowk$
makes sense for any ground field $k$.
Note that for $k$ finite one can still read off the Hasse--Weil zeta function
from it.

As pointed out by Andr\'e in Section 4.3 of \cite{Andre} and Chapter 13 of
\cite{AndreBuch}, 
varieties with a finite dimensional Chow motive in the sense of Kimura
\cite{Kimura} and
O'Sullivan (i.e. whose Chow motive is the sum of two Chow motives $X^+$  
and $X^-$ such
that $\Alt^{i}(X^+)=0$ for $i\gg 0$ and $\Sym^{i}(X^-)=0$ for $i\gg 0$) have a
rational zeta function with coefficients in $\Chowk$. More precisely,
as $\Alt^i$ is the opposite structure to $\Sym^i$ (compare Section
\ref{lambda}), 
\begin{equation}\label{rational}
Z_X(T)=\frac{P(T)}{Q(-T)}\text{ in }\Chowk[[T]], 
\end{equation}
where $P(T)=\sum_{i\ge 0}[\Sym^i(X^-)]T^i$ and 
$Q(T)=\sum_{i\ge 0}[\Alt^i(X^+)]T^i$ are polynomials and moreover $Q(T)$ is
invertible in $\Chowk[[T]]$. For example, this holds 
for an abelian variety over $k$. 

In Chapter 13 of \cite{AndreBuch}, Andr\'e writes: \selectlanguage{french} \flqq Nous laissons au lecteur le plaisir de sp\'eculer sur 
d'\'eventuelles \'equations fonctionelles\dots \frqq
\selectlanguage{english}

\medskip
 
In this note, we consider functional equations for zeta functions with
coefficients in $\Chowk$, where $k$ is an arbitrary field.
Using the well known decomposition of the Chow motive of an
abelian variety, we prove
\begin{prop}[Proposition \ref{abelsch}]
Let $A$ be an abelian variety of dimension $g$ over $k$, and denote by
$Z_A(T)=\sum_{i=0}^{\infty} [\Sym^ih(X)]T^i\in\Chowk[[T]]$ its zeta function
with values in $\Chowk$. 
Then $$Z_A(\frac{1}{\tate^g T})=Z_A(T).$$ More
precisely, 
$Z_A(T)$ can be written as $Z_A(T)=\frac{P^A(T)}{Q^A(-T)}$ as in
Equation \ref{rational}
in such a way that
$P^A(T),Q^A(T)\in1+T\Chowk[T]$ satisfy the expected functional equations
$$P^A(\frac{1}{\tate^g T})=T^{-f}\tate^{-\frac{gf}{2}}P^A(T)\text{ and }
Q^A(\frac{1}{\tate^g T})=T^{-e}\tate^{-\frac{ge}{2}}Q^A(T)$$
in $\Chowk[T,T^{-1}]$, where $e=f=2^{2g-1}$.
\end{prop} 
Furthermore, in Proposition \ref{produkte}, we show that 
having a rational zeta function satisfying a 
functional equation is
preserved by taking products. 

To this end, in Section \ref{lambdaknull}, we investigate the
$\lambda$--structure on $\Chowk$.

\subsubsection*{Acknowledgments} I am indebted to H\'el\`ene Esnault 
for her invaluable suggestions and comments and constant
encouragement. I thank Niko Naumann
for interesting discussions. I am grateful for the support of the 
\selectlanguage{german} DFG--Schwerpunkt \glqq Globale Methoden in der Komplexen Geometrie\grqq. 
\selectlanguage{english} 

\section{$\lambda$--Rings}
\label{lambda}
Recall the notion of a $\lambda$--ring. For more details,
see for example Chapter I of Atiyah and Tall, \cite{AtiyahTall}.

A ring $A$ endowed with operations $\lambda^r$ for $r\in\nat$ 
such that $\lambda^0(a)=1$, $\lambda^1(a)=a$ and 
$\lambda^r(a+b)=\sum_{i+j=r} \lambda^i(a)\lambda^j(b)$ is called a
\emph{$\lambda$--ring}. 
This is equivalent to the datum of a group homomorphism $\lambda_t:
(A,+)\longrightarrow (1+tA[[t]],\cdot), a\mapsto 1+\sum_{r\ge 1}
\lambda^r(a)t^r$ such that $\lambda^1(a)=a$.

The \emph{opposite} $\lambda$--structure on $A$ is given by 
$\sigma_t(a)=(1+\sum_{r\ge 1}
\lambda^r(a)(-t)^r)^{-1}$. Explicitely, $\sigma^r$ is given recursively by 
$$\sigma^r(a)-\sigma^{r-1}(a)\lambda(a)+\dots+(-1)^r\lambda^r(a)=0
\text{ for }r\ge 1.$$

$B=1+tA[[t]]$ itself carries the structure of a $\lambda$--ring: 

Denote by
$\sigma_i^N$ the elementary symmetric polynomials in $\xi_1,\dots,\xi_N$ 
and by $s_i^N$ the
elementary symmetric polynomials in $x_1,\dots,x_N$.
Let $P_n(\sigma_1^N,\dots,\sigma_n^N,s_1^N,\dots,s_n^N)$ be
the coefficient of $t^n$ in $\prod_{1\le i,j\le N}(1+\xi_i x_j t)$, where
$N\ge n$,  
and $P_{n,r}(\sigma_1^N,\dots,\sigma_{rn}^N)$ 
the coefficient of $t^n$ in
$\prod_{1\le i_1<\dots<i_r\le N}(1+\xi_{i_1}\dotsm\xi_{i_r}t)$, where $N\ge rn$.

Addition on $B$ is given by multiplication,
multiplication $\circ$ is given by 
$$(1+\sum_{k\ge 1} a_kt^k) \circ 
(1+\sum_{l\ge 1} b_lt^l) = 1+\sum_{n\ge 1}P_n(a_1,\dots,a_n;b_1,\dots,b_n)t^n$$
with neutral element $1+t$, 
and the $\lambda$--structure is given by
$$\Lambda^r(1+\sum_{k\ge 1} a_kt^k )=1+\sum_{n\ge 1}
P_{n,r}(a_1,\dots,a_{rn})t^n.$$
The $\lambda$--ring $A$ is called \emph{special}, if $\lambda_t$ is a
homomorphism of $\lambda$--rings.

\begin{rem}
The $\lambda$--structure on $B$ may be given in a more sophisticated manner
without writing down the universal polynomials $P_n$ and $P_{n,r}$ explicitly, 
compare Section I.1 of \cite{AtiyahTall}. But we will need the precise shape
of $P_n$ and $P_{n,r}$ in Sections \ref{abschnabelsch} and \ref{abschnprodukte}.
\end{rem}

\begin{rem}\label{erz}
A group homomorphism $\varphi:A\longrightarrow B$ between
$\lambda$--rings is a homomorphism of $\lambda$--rings if there is a set of
group generators $S\subseteq A$ such that
$\varphi(ab)=\varphi(a)\varphi(b)$ for all $a,b\in S$ and for all $r\in \nat$
and $a\in S$ we have $\varphi(\lambda^r(a))=\lambda^r(\varphi(a))$. Compare
\cite{LarsenLunts2}, Lemma 4.4.
\end{rem}

\section{Curves}
As a motivation, let us briefly review the situation for curves.

First, we consider the zeta function associated to the universal Euler
characteristic with compact support.
Let $C$ be a geometrically irreducible
smooth projective curve of genus $g$ over a field $k$. Denote its $i$-th
symmetric power by $\Sym^i(C)$.
The zeta function of $C$ is defined as

$$\zeta_C(T) =\sum_{i=0}^{\infty} [\Sym^i(C)]T^i\text{ in }\Knull(\Vark)[[T]].$$
Here $\Knull(\Vark)$ is 
the value group of the universal Euler characteristic 
with compact support, i.e. the free abelian group on isomorphism classes of
varieties over $k$ modulo the relations $[X]=[X-Y]+[Y]$, where $Y\subset X$
closed. It carries a commutative ring structure 
induced by the product of varieties. By abuse of notation, we denote the
class of the affine line by $\tate$. A stratification argument shows that we
get the same Grothendieck ring if we take classes of \emph{quasi-projective}
varieties. 

If there is a line bundle of degree $1$ on $C$, Kapranov shows that
 $(1-T)(1-\tate T)\zeta_C(T)$ is a polynomial of degree $2g$, 
and that the zeta function satisfies the functional equation
$$\zeta_C(\frac{1}{\tate T})=\tate^{1-g}T^{2-2g}\zeta_C(T)$$
in $\cM_k((T))$, where
$\cM_k:= \Knull(\Vark)[\tate^{-1}]$.

Let us give a slight reformulation of Kapranov's argument.

As pointed out by Larsen and Lunts in \cite{LarsenLunts2}, 
the symmetric powers $\Sym^i$ of quasi-projective varieties
induce the structure 
of a $\lambda$--ring on $\Knull(\Vark)$ and, since $\Sym^i(\tate[X])
=\tate^i\Sym^i([X])$ (see G\"ottsche, \cite{Goettsche}), also on $\cM_k$ via 
$\Sym^i(\tate^k[X]):=\tate^{ik}\Sym^i([X])$. In these terms, 
$\zeta_X(T)=\lambda_t([X])$. Since $\lambda_t$ is a group homomorphism, we get
$$\zeta_C(T)=
(\sum_{i=0}^{\infty}\Sym^i([C]-[\proj{1}])T^i) \zeta_{\proj{1}}(T).$$
As $\zeta_{\proj{1}}(T)=\zeta_1(T)\zeta_{\tate}(T) 
= \frac{1}{1-T} \frac{1}{1-\tate T}$, multiplying this equation by $(1-T)(1-\tate T)$ yields
$\Sym^i([C]-[\proj{1}]) = [\Sym^i(C)]-[\Sym^{i-1}(C)][\proj{1}] +
[\Sym^{i-2}(C)]\tate$ for $i\ge 2$.

If $i>2g$, this expression vanishes, because
for $j>2g-2$ the morphism $\Sym^j(C)\longrightarrow \Pic^j(C) \cong \Pic^0(C)$ 
is a Zariski fibration with fiber $\proj{j-g}$ (we still assume that there is
a line bundle of degree $1$ on $C$). Therefore, 
 $(1-T)(1-\tate T)\zeta_C(T)$ is a polynomial of degree $2g$.

For the functional equation we need to show for $g\le i\le 2g$ that
$\Sym^i([C]-[\proj{1}])=\tate^{i-g}\Sym^{2g-i}([C]-[\proj{1}])$:
Consider the morphism $\Sym^j(C)\longrightarrow \Pic^j(C)$. It is a piecewise 
Zariski fibration with fiber $\Proj{H^0(L)}$ over $L$.
There is an isomorphism $\Pic^i(C)\cong\Pic^{2g-2-i}(C)$, 
$L\mapsto\omega_C\otimes L^\vee$.
By Riemann--Roch, $h^0(L)-h^0(\omega_C\otimes L^\vee)=\deg L+1-g$.
Therefore, 
$$[\Sym^i(C)]-\tate^{i-g+1}[\Sym^{2g-2-i}(C)]=[\proj{i-g}][\Pic^i(C)]\text{
for }g\le i\le 2g-2.$$
Using $\Pic^j(C) \cong \Pic^0(C)$ again and adding up we conclude
$\Sym^i([C]-[\proj{1}])=\tate^{i-g}\Sym^{2g-i}([C]-[\proj{1}])$.

Actually, the equation 
$\zeta_{\proj{1}}(T) = \frac{1}{1-T} \frac{1}{1-\tate T}$ can be rephrased by
saying that $\Alt^i(1)=\Alt^i(\tate)=0$ for $i\ge 2$, where 
$\Alt^i$ is the opposite $\lambda$--structure on $\Knull(\Vark)$.

Hence, $[C]$ can be written as the sum of two terms $x^++x^-$, 
where $\Alt^{i}(x^+)=0$ for $i\gg 0$ and $\Sym^{i}(x^-)=0$ for $i\gg 0$ and
furthermore, $P^{x^-}(T)=\sum_i\Sym^i(x^-)T^i$ and
$Q^{x^+}(T)=\sum_i\Alt^i(x^+)T^i$ satisfy the expected functional equations.

\begin{rem}
If $C$ carries a line bundle of degree $d$, a similar calculation shows that
$(1-T^d)(1-\tate^dT^d)\zeta_C(T)$ is a polynomial and that the functional
equation still holds, as pointed out by Kapranov.
\end{rem}

Now let us consider the zeta function of $C$ with values in $\Chowk$, 
$$Z_C(T)=\sum_{i=0}^{\infty} [\Sym^i h(C)]T^i.$$
The Chow motive of $C$ has a decomposition
$$h(C)=\eins\oplus h^1(C)\oplus\tate,$$ where
$\Sym^{i}h^1(C)=0$ for $i>2g$ and due to a result of K\"unnemann in \cite{kuennemann}, 

$\Sym^ih^1(C)\cong\tate^{g-i}\otimes\Sym^{2g-i}h^1(C)$ (compare Section 3.3 in
\cite{BanoDiss}), 
therefore the zeta function of $C$ with values in $\Chowk$ 
is also rational and
satisfies the expected functional equation. In characteristic zero, this
follows from the properties of the zeta function with values in $\cM_k$ 
(compare Section \ref{lambdachow}). 

\section{$\lambda$--Structures on the Grothendieck ring of Chow motives}

\label{lambdaknull}
For the rest of the paper, we will restrict ourselves to the study of 
zeta functions with values in $\Chowk$. In fact, the properties of the
$\lambda$--structure on $\Chowk$ which we need hold for the Grothendieck 
ring of any (pseudo-abelian) $\rat$-linear tensor category.

\subsection{Schur functors}

Let us recall some facts from Deligne, \cite{Deligne}, Section 1.

Let $\kappa$ be a field of characteristic zero, let 
$\cA$ be a \emph{$\kappa$-linear tensor category}, i.e. a symmetric monoidal category,
which is additive, pseudo-abelian and $\kappa$-linear such that $\otimes$ is
$\kappa$-bilinear. 

If $V$ is a finite dimensional $\kappa$-vector space and $X$ is an object of $\cA$,
there are objects $V\otimes X$ and $\cHom(V,X)$ of $\cA$ natural in $V$ and $X$ 
such that
$$\Hom(V\otimes X, Y)=\Hom(V,\Hom(X,Y))$$ and
$$\Hom(Y,\cHom(V,X))=\Hom(V\otimes Y,X).$$

There is a natural isomorphism $\cHom(V,X)\cong V^\vee\otimes X$. The choice
of a basis of $V$ yields (non-canonical) isomorphisms
$\cHom(V,X) \cong X^{\oplus \dim V} \cong V \otimes X$.

If a finite group $G$ acts on $X$, we define $X^G$ as the image of the
projector $\frac{1}{|G|}\sum_{g\in G}g\in\End(X)$.

If $G$ acts on $V$ and on $X$, it acts on $\cHom(V,X)$ and we define
$\cHom_G(V,X)$ as $\cHom(V,X)^G$. Note that 
$\Hom(Y,\cHom_G(V,X))=\Hom_G(V,\Hom(Y,X))$.

If all irreducible representations of $G$ over $\overline{\kappa}$ 
are already defined over $\kappa$ we have 
$\kappa[G]\cong\prod \End_\kappa(V_\lambda),$
where $V_\lambda$ runs through a system of representatives for the isomorphism
classes of irreducible representations.
Therefore, 
\begin{equation}\label{zerlegung}
X\cong\cHom_G(\kappa[G],X)\cong\bigoplus 
V_\lambda\otimes\cHom_G(V_\lambda,X),
\end{equation}
where the $G$-action on $X$ corresponds to the $G$-action on the outer
$V_\lambda$ on the right hand side.

There is a natural isomorphism 
$$\cHom_{G\times H}(V\otimes W, X\otimes Y)\cong
\cHom_G(V,X)\otimes\cHom_H(W,Y).$$ 
Under this isomorphism, the $S_n$-action on $\cHom_G(V,X)^{\otimes n}$
corresponds to the $S_n$-action on $\cHom_{G^n}(V^{\otimes n},X^{\otimes n})$
induced by the actions on $V^{\otimes n}$ and $X^{\otimes n}$.

Furthermore, if we have a short exact sequence of finite groups
$$1\longrightarrow K\longrightarrow G\longrightarrow H\longrightarrow 1,$$
and $V$ is a representation of $H$ while $W$ is a representation of $G$, which 
also acts on $X\in\cA$,
we have a natural isomorphism 
$$\cHom_G(V\otimes W,X)\cong\cHom_H(V,\cHom_K(W,X)).$$

Finally, if $G$ acts on $X$ and $V$ is a representation of a subgroup $H<G$,
$$\cHom_H(V,X)\cong\cHom_G(\Ind^G_H V,X).$$

If $G$ is the symmetric group $S_n$ and $V_\lambda$ is an irreducible
representation of $S_n$, indexed by a partition $\lambda$ of $n=|\lambda|$,
$$S_\lambda(X):=\cHom_{S_n}(V_\lambda,X^{\otimes n})$$
is called \emph{Schur functor}.
For the trivial representation $\Triv(S_n)$ we get 
$$\Sym^n(X):=S_{(n)}(X)=\im(\frac{1}{n!}\sum \sigma)\subseteq X^{\otimes n},$$ 
for the alternating representation $\Sign(S_n)$ we obtain 
$$\Alt^n(X):=S_{(1^n)}(X)=\im(\frac{1}{n!}\sum (-1)^\sigma \sigma)\subseteq X^{\otimes n}.$$ 

By \ref{zerlegung}, there is a canonical isomorphism
$$X^{\otimes n}\cong\bigoplus_{|\lambda|=n} V_\lambda\otimes S_\lambda(X),$$
where the  $S_n$--action on $X^{\otimes n}$ corresponds to the action on
$V_\lambda$ on the right hand side. Note that in particular
\begin{equation}\label{schurverschwind}
S_\lambda(X)=0\text{ for }\lambda\ne (n), 
\end{equation}
if $\Sym^n(X)= X^{\otimes n}$.

If $V_{\mu_i}$, $i=1,\dots, r$ are irreducible representations of $S_{n_i}$ and
$V_\lambda$ is an irreducible representation of $S_n$, where $n=\sum
n_i$, we denote by $$[\lambda:\mu_1,\dots,\mu_r]$$ the multiplicity of
$\bigotimes V_{\lambda_i}$ in 
$\Rest^{S_n}_{\prod S_{n_i}} V_\lambda$ (which equals the multiplicity of 
$V_\lambda$ in $\Ind^{S_n}_{\prod S_{n_i}}\bigotimes V_{\lambda_i}$
by Frobenius reciprocity).

With this notation, we get
\begin{align}
S_\mu(X)\otimes S_\nu(X)&\cong \cHom_{S_m\times S_n}(V_\mu\otimes
V_\nu,X^{\otimes m+n})\notag\\ &\cong \cHom_{S_{m+n}}(\Ind^{S_{m+n}}_{S_m\times
S_n}V_\mu\otimes V_\nu,X^{\otimes m+n})\label{malinnen}\\
&\cong\bigoplus_{|\lambda|=|\mu|+|\nu|}[\lambda:\mu,\nu]S_\lambda(X)\label{malaussen}\\
S_\lambda(X\oplus Y)&\cong\bigoplus_{|\mu|+|\nu|=|\lambda|}[\lambda:\mu,\nu]
S_\mu(X)\otimes S_\nu(Y)\label{add}\\
S_\lambda(X\otimes Y)&\cong\bigoplus_{|\mu|=|\nu|=|\lambda|}[V_\mu\otimes
  V_\nu:V_\lambda]S_\mu(X)\otimes S_\nu(Y)\label{schurprodukt}
\end{align}

Furthermore, for a $S_m$--representation $V$ and a $G$--representation $W$, 
\begin{equation}\label{kranz}
\cHom_{S_m}(V,\cHom_{G^m}(W^{\otimes m},X^{\otimes m}))
\cong \cHom_{S_m\ltimes G^m}(V\otimes{W}^{\otimes m}, X^{\otimes m}).
\end{equation}

In particular, for $|\mu|=m$ and $|\nu|=n$, 
\begin{align*}
S_\mu(S_\nu(X))&=
 \cHom_{S_m\ltimes{S_n}^m}(V_\mu\otimes{V_\nu}^{\otimes m}, X^{\otimes nm})\\
&\cong \cHom_{S_{nm}}(\Ind^{S_{nm}}_{S_m\ltimes {S_n}^ m}
V_\mu\otimes{V_\nu}^{\otimes m}, X^{\otimes nm}).
\end{align*}

\subsection{The $\lambda$--structure}

Denote by $\Knull(\cA)$ the free abelian group on isomorphism classes 
$[X]$ of objects of $\cA$ modulo the relations $[X\oplus Y]=[X]+[Y]$. It is
the Grothendieck group associated to the abelian monoid of isomorphism classes 
of objects in $\cA$ with direct sum.
The
tensor product of $\cA$ induces a commutative 
ring structure on $\Knull(\cA)$. We call $\Knull(\cA)$ the \emph{Grothendieck
ring of $\cA$}.

Note that for any $X\in \cA$ with $G$--action we obtain a group homomorphism 
from the Grothendieck group of $G$--representations to $\Knull(\cA)$ sending
a representation $V$ to $\cHom_G(V,X)$.

\begin{lem}
The exterior powers $\Alt^n$ induce a special $\lambda$--ring structure on
$\Knull(\cA)$ with opposite $\lambda$--structure given by the symmetric powers
$\Sym^n$.
\end{lem}

\begin{proof}
Due to Equation \ref{add} and the Littlewood--Richardson rule (see e.g.
\cite{FultonHarris}, Appendix A),
$$[X]\mapsto 1+\sum_{n\ge 1}[\Alt^n(X)]t^n$$ induces a 
$\lambda$--ring 
structure on $\Knull(\cA)$.
The fact that the opposite structure is given by $\Sym^n$ follows from
Equation \ref{malaussen} and the Littlewood--Richardson rule
or more precisely from the fact that for $i,j\ge 1$ we have
$$\Sym^i(X)\otimes\Alt^j(X)\cong S_{(i+1,1^{j-1})}(X)\oplus 
S_{(i,1^j)}(X).$$

To show that the $\lambda$--structure given by $\Alt^i$ is special, 
we use an argument by Larsen and Lunts from \cite{LarsenLunts2},
Theorem 5.1, in a
slightly more general setting.

Recall one possible description of the free special $\lambda$-ring $R$ on one generator:
$R=\bigoplus_{n\ge 0} R_n$, where $R_n$ is the  
representation ring over $\kappa$ of the
symmetric group $S_n$ (with the convention that $S_0$ is the trivial group and
$R_0$ therefore is $\ganz$).
It has a $\ganz$-basis consisting of the elements $(n,V_\nu)$, where $V_\nu$ is an irreducible
$S_n$-representation. The product is given by
$$(m,V_\mu)(n,V_\nu)=(m+n,\Ind^{S_{m+n}}_{S_m\times S_n}V_\mu\otimes V_\nu),$$
while the $\lambda$-structure is is given by
$$\lambda^r(n,V_\nu)=(rn,\Ind^{S_{rn}}_{S_r\ltimes {S_n}^r}\Sign(S_r)\otimes 
{V_\nu}^{\otimes r}).$$
Its generator as a $\lambda$-ring is $(1,\kappa)$.

$R\otimes_{\ganz}R$ is the free special $\lambda$-ring on two generators.
It has a $\ganz$-basis consisting of the elements $(n_1,n_2,V_{\nu_1}\otimes V_{\nu_2})$, where 
$V_{\nu_i}$ is an irreducible
$S_{n_i}$-representation.
The product is given by
\begin{multline*}
(m_1,m_2,V_{\mu_1}\otimes V_{\mu_2})(n_1,n_2,V_{\nu_1}\otimes V_{\nu_2})=\\
(m_1+n_1,m_2+n_2,\Ind^{S_{m_1+n_1}\times S_{m_2+n_2}}_{(S_{m_1}\times S_{n_1})\times(S_{m_2}\times S_{n_2})}(V_{\mu_1}\otimes 
V_{\nu_1})\otimes(V_{\mu_2}\otimes 
V_{\nu_2})),
\end{multline*}
while the $\lambda$-structure is is given by
\begin{multline*}
\lambda^r(n_1,n_2,V_{\nu_1}\otimes V_{\nu_2})=\\
(rn_1,rn_2,\Ind^{S_{rn_1}\times S_{rn_2}}_{S_r\ltimes ({S_{n_1}}^r\times {S_{n_2}}^r)}\Sign(S_r)\otimes 
{V_{\nu_1}}^{\otimes r}\otimes {V_{\nu_2}}^{\otimes r}).
\end{multline*}

Now let $X_1,X_2$ be two objects of $\cA$. Then, by \ref{malinnen} and
\ref{kranz},
$$(n_1,n_2,V_{\nu_1}\otimes V_{\nu_2})\mapsto S_{\nu_1}(X_1)\otimes
S_{\nu_2}(X_2)$$ defines a $\lambda$-ring homomorphism $R\otimes_\ganz
R\longrightarrow
\Knull(\cA)$, hence every pair $[X_1],[X_2]$ is contained in a special
$\lambda$-subring of $\Knull(\cA)$.

Therefore, $\lambda_t:\Knull(\cA)\longrightarrow 1+t\Knull(\cA)[[t]]$
satisfies $$\lambda_t(xy)=\lambda_t(x)\circ\lambda_t(y)$$ and
$$\lambda_t(\lambda^r x)=\Lambda^r(\lambda_t(x))$$
for elements $x=[X]$ and $y=[Y]$ and due to Remark \ref{erz} therefore for all $x,y\in\Knull(\cA)$.
\end{proof}

We will need some more identities relating symmetric and exterior powers.

For a representation $V$ of $S_n$, let $V':=\Sign(S_n)\otimes V$.
Note that 
\begin{equation*}
\Sign(S_{mn})\otimes\Ind^{S_{mn}}_{S_m\ltimes{S_n}^m}
V\otimes W ^{\otimes m}\cong
\begin{cases}
\Ind^{S_{mn}}_{S_m\ltimes{S_n}^m}
V'\otimes W'^{\otimes m}& \text{if }n\text{ is odd}\\
\Ind^{S_{mn}}_{S_m\ltimes{S_n}^m}
V\otimes W'^{\otimes m}& \text{if }n\text{ is even,}
\end{cases}
\end{equation*}
because 
\begin{equation*}
(-1)^{\iota(\sigma)}=
\begin{cases}
(-1)^\sigma & \text{if }n\text{ is odd}\\
1 & \text{if }n\text{ is even,}
\end{cases}
\end{equation*}
where $\iota:S_m\hookrightarrow S_{mn}$.
Furthermore,
$$\Sign(S_{m+n})\otimes\Ind^{S_{m+n}}_{S_m\times S_n}V\otimes W
\cong \Ind^{S_{m+n}}_{S_m\times S_n}V'\otimes W'.$$
Therefore, for $n$ odd, we obtain the following equation in $R$:
\begin{equation*}
\begin{split}
\Ind^{S_{mn}}_{S_m\ltimes{S_n}^m}
\Triv(S_m)\otimes \Triv(S_n)^{\otimes m}
&=\Ind^{S_{mn}}_{S_m\ltimes{S_n}^m}
\Sign(S_m)'\otimes \Sign(S_n)'^{\otimes m}\\
&=\Sign(S_{mn})\otimes P_{m,n}(\Sign(S_1),\dots,\Sign(S_{mn}))
\\
&=  P_{m,n}(\Triv(S_1),\dots,\Triv(S_{mn})).
\end{split}
\end{equation*}
Similarly, for $n$ even, we get

$$\Ind^{S_{mn}}_{S_m\ltimes{S_n}^m}\Sign(S_m)\otimes \Triv(S_n)^{\otimes m}
= P_{m,n}(\Triv(S_1),\dots,\Triv(S_{mn})).$$

Hence for every $x=[X]\in\Knull(\cA)$ we get
\begin{eqnarray}
\Sym^m(\Sym^n(x))&=P_{m,n}(\Sym^1(x),\dots,\Sym^{mn}(x))&\text{if }n\text{ is odd,}\label{ungerade}\\
\Alt^m(\Sym^n(x))&=P_{m,n}(\Sym^1(x),\dots,\Sym^{mn}(x))&\text{if }n\text{ is even.}\label{gerade}
\end{eqnarray}

Now consider the two generators $e_1=(1,0,\kappa\otimes\kappa)$ and
$e_2=(0,1,\kappa\otimes\kappa)$ 
of $R\otimes_\ganz R$.

We know that

$$\lambda^n(e_1e_2)=P_n(\lambda^1(e_1),\dots,\lambda^n(e_1),\lambda^1(e_2),\dots,
\lambda^n(e_2)),$$
where $\lambda^i(e_1)=(i,0,\Sign(S_i)\otimes\kappa)$ and
$\lambda^i(e_2)=(0,i,\kappa\otimes\Sign(S_i))$.
On the other hand,
$$\lambda^n(e_1e_2)=\sum_{|\mu|=|\nu|=n}[V_\mu\otimes V_\nu:\Sign(S_n)](n,n,V_\mu \otimes V_\nu).$$
As $[V_\mu\otimes V_\nu:\Sign(S_n)]=[V_\mu'\otimes V_\nu':\Sign(S_n)]$, 
applying Equations \ref{malaussen} and \ref{malinnen}, 
for $x=[X],y=[Y]\in\Knull(\cA)$ we get
\begin{equation}
\Alt^n(xy)=
P_n(\Sym^1(x),\dots,\Sym^n(x),\Sym^1(y),\dots,
\Sym^n(y)).\label{symsym}
\end{equation}
Similarly, as $[V_\mu\otimes V_\nu:\Triv(S_n)]=[V_\mu'\otimes
V_\nu:\Sign(S_n)]=[V_\mu\otimes V_\nu':\Sign(S_n)]$, we have
\begin{align}
\Sym^n(xy)&=
P_n(\Sym^1(x),\dots,\Sym^n(x),\Alt^1(y),\dots,
\Alt^n(y))\label{symalt}\\
&=P_n(\Alt^1(x),\dots,\Alt^n(x),\Sym^1(y),\dots,
\Sym^n(y)).\label{altsym}
\end{align}

\subsection{The Grothendieck ring of Chow motives}
\label{lambdachow}
Everything in this section applies to the $\rat$--linear tensor category
$\mathit{CM}_k$ of
Chow motives over a field $k$ with rational coefficients as in Manin,
\cite{Manin} or Scholl, \cite{Scholl}, where the equivalence relation on
cycles is rational equivalence. 
Note that $\Sym^n \tate = \tate^{\otimes n}$ and therefore $S_\lambda(\tate)=0$ for 
$\lambda\ne (n)$ due to \ref{schurverschwind}.
Hence it follows from \ref{schurprodukt} that
\begin{equation}
S_\lambda(\tate \otimes M)\cong\tate^{\otimes|\lambda|}\otimes S_\lambda(X).
\label{schurtate}
\end{equation}
We denote the Grothendieck ring
of $\mathit{CM}_k$ by $\Chowk$.
If $k$ is of characteristic zero, Gillet and Soul\'e as a corollary from
\cite{GilletSoule} 
get a ring homomorphism from $\cM_k$ to $\Chowk$
such that for a smooth projective variety $X$ the class $[X]$ of $X$ is sent
to $[h(X)]$, where $h(X)$ is 
the Chow motive of $X$.
Del Ba{\~n}o and Navarro Aznar have shown in \cite{BanoNavarro} that
for a finite group $G$ acting on $X$, the class $[X/G]$ of the quotient is 
sent to $[h(X)^G]$, where $h(X)^G$ is the image of the projector
$\frac{1}{|G|}\sum_{g\in G}g$ in $h(X)$. 
In particular, the ring homomorphism $\cM_k\longrightarrow\Chowk$
is actually a homomorphism of
$\lambda$--rings. 

\section{Abelian varieties}
\label{abschnabelsch}
The aim of this section is to prove the following
\begin{prop} \label{abelsch}
Let $A$ be an abelian variety of dimension $g$ over $k$, and denote by
$Z_A(T)=\sum_{i=0}^{\infty} [\Sym^ih(X)]T^i\in\Chowk[[T]]$ its zeta function
with values in $\Chowk$. Then $Z_A(\frac{1}{\tate^g T})=Z_A(T)$. More
precisely, $Z_A(T)$ can be written as $Z_A(T)=\frac{P^A(T)}{Q^A(-T)}$, 
such that
$P^A(T),Q^A(T)\in1+T\Chowk[T]$ satisfy the expected functional equations
$$P^A(\frac{1}{\tate^g T})=T^{-f}\tate^{-\frac{gf}{2}}P^A(T)\text{ and }
Q^A(\frac{1}{\tate^g T})=T^{-e}\tate^{-\frac{ge}{2}}Q^A(T)$$
in $\Chowk[T,T^{-1}]$, where $e=f=2^{2g-1}$.
\end{prop}

\begin{proof}
As shown by Beauville in \cite{Beauville} and Deninger and Murre in
\cite{DeningerMurre}, the Chow motive of an abelian variety of dimension $g$ 
is canonically isomorphic to a sum
$$h(A)\cong\bigoplus_{0\le i\le 2g}h^i(A),$$
where $h^i(A)\cong\Sym^i(h^1(A))$ (in particular, $\Sym^i(h^1(A))=0$ for
$i>2g$), and multiplication by $n$ acts on $h^i(A)$ as $n^i$. Furthermore,
$h^0(A)=1$ and $h^{2g}(A)=\tate^g$.
Therefore, the zeta function of $A$ with values in $\Chowk$ equals
$$Z_A(T)
=\frac{\prod_{\substack{0\le n\le 2g\\n\text{
	odd}}}P_n^A(T)}{\prod_{\substack{0\le n\le 2g\\n\text{
	even}}}Q_n^A(-T)},$$
where
\begin{align*}
P^A_n(T)&:=\sum_{m\ge 0}[\Sym^m(\Sym^n(h^1(A)))]T^m,\\
Q^A_n(T)&:=\sum_{m\ge 0}[\Alt^m(\Sym^n(h^1(A)))]T^m.
\end{align*}
Denote by $\sigma_1^N,\dots,\sigma_N^N$ 
the elementary symmetric functions in $\xi_1,\dots,\xi_N$.
From the commutativity of the diagram
$$\xymatrix{\ganz[\sigma_1^K,\dots,\sigma_K^K]\ar[r]\ar@{^{(}->}[d]
&\ganz[\sigma_1^k,\dots,\sigma_k^k]\ar@{^{(}->}[d]\\
\ganz[\xi_1,\dots,\xi_K]\ar[r]
&\ganz[\xi_1,\dots,\xi_k]}$$
where $k\le K$,  
$\sigma_l^K\mapsto 0$ for $k<l\le K$ and to $\sigma_l^k$ for $l\le k$ and
$\xi_l\mapsto 0$ for $k<l\le K$ under the horizontal maps,
it follows that
$$q^g_n(t):=\sum_{m\ge 0}P_{m,n}(\sigma_1^{2g},\dots,\sigma_{2g}^{2g},0,\dots,0)t^m=\prod_{1\le
i_1<\dots<i_n\le 2g}(1+\xi_{i_1}\dotsm\xi_{i_n} t)$$
is a polynomial of degree $b^g_n=\binom{2g}{n}$. For convenience, let us
denote $\sigma^{2g}_{2g}$ by $\sigma$.
The polynomial $q^g_n(t)$ satisfies 
\begin{align*}
q^g_n(\frac{1}{\sigma t})&=\prod_{1\le i_1<\dots<i_n\le
2g}(1+\xi_{i_1}\dotsm\xi_{i_n} \frac{1}{\sigma t})\\
&=(\frac{1}{\sigma t})^{b^g_n}\prod_{1\le i_1<\dots<i_n\le
2g}\xi_{i_1}\dotsm\xi_{i_n}\prod_{1\le i_1<\dots<i_n\le
2g}(1+\frac{\sigma}{\xi_{i_1}\dotsm\xi_{i_n}}t)\\ &=
(\frac{1}{\sigma t})^{b^g_n}{\sigma}^{\frac{b^g_n n}{2g}}q^g_{2g-n}(t).
\end{align*}
Therefore, for $n$ odd, it follows from \ref{ungerade} that 
$P^A_n(T)$
is a polynomial of degree
$b_n^g$ and 
$$P^A_n(\frac{1}{\tate^g
T})=(\frac{1}{\tate^g T})^{b^g_n}\tate^{\frac{b^g_n
n}{2}}P^A_{2g-n}(T).$$ In particular, 
$P^A(T):=\prod_{\substack{0\le n\le 2g\\n\text{
	odd}}}P_n^A(T)$ satisfies
$P^A(\frac{1}{\tate^g T})=T^{-f}\tate^{-\frac{gf}{2}}P^A(T).$

On the other hand, for $n$ even, we deduce from \ref{gerade} that
$Q^A_n(T)$
is a polynomial of degree
$b_n^g=\binom{2g}{n}$ and satisfies
$$Q^A_n(\frac{1}{\tate^g
T})=(\frac{1}{\tate^g T})^{b^g_n}\tate^{\frac{b^g_n
n}{2}}Q^A_{2g-n}(T),$$
hence $Q^A(T):=\prod_{\substack{0\le n\le 2g\\n\text{
	even}}}Q_n^A(T)$ satisfies
$Q^A(\frac{1}{\tate^g T})=T^{-e}\tate^{-\frac{ge}{2}}Q^A(T).$
\end{proof}

\begin{rem}
An easy calculation using Equation \ref{schurtate} and the decomposition of
the motive of a blow--up as e.g. in \cite{Manin}, Section 9,
shows that the property of having a rational zeta function
satisfying a functional equation is closed under blow--ups along smooth 
centers satisfying a functional equation.

More precisely, suppose that $X$ is an $n$--dimensional smooth projective
variety such that $[h(X)]=[X^+]+[X^-]$, where $[\Alt^i(X^+)]=0$ for $i>e(X^+)$
and $[\Sym^i(X^-)]=0$ for $i>f(X^-)$.
Let $Q^{X^+}(T)=\sum_{i\ge 0} [\Alt^i(X^+)]T^i$ and 
$P^{X-}(T)=\sum_{i\ge 0} [\Alt^i(X^-)]T^i$.
Suppose furthermore that
$$Q^{X^+}(\frac{1}{\tate^n T})=T^{-e(X^+)}\tate^{-\frac{ne(X^+)}{2}}Q^{X^+}(T)
\text{ in }\Chowk[T,T^{-1}]$$
and 
$$P^{X^-}(\frac{1}{\tate^n
T})=T^{-f(X^-)}\tate^{-\frac{nf(X^-)}{2}}P^{X^-}(T)\text{ in
}\Chowk[T,T^{-1}],$$ and likewise for a smooth closed subvariety $Y$ of $X$ of
pure codimension $d$.  
Then the same holds for the blow--up $\Bl{Y}{X}$ of $X$ along
$Y$, where $(\Bl{Y}{X})^+=X^+\oplus\bigoplus_{i=1}^{d-1}\tate^i\otimes Y^+$,
$(\Bl{Y}{X})^-=X^-\oplus\bigoplus_{i=1}^{d-1}\tate^i\otimes Y^-$,
$e((\Bl{Y}{X})^+)=e(X^+)+(d-1)e(Y^+)$ and
$f((\Bl{Y}{X})^-)=f(X^-)+(d-1)f(Y^-)$.
\end{rem}
\begin{rem}
For Kummer surfaces $X$, an explicit calculation of $[h(X)]$ 
(we know how multiplication
by $-1$ acts on the Chow motive of an abelian variety) yields 
$[\Alt^i(h(X))]=0$
for $i>24$ and the expected functional equation
$Q^{h(X)}(\frac{1}{\tate^2 T})=T^{-24}\tate^{-24}Q^{h(X)}(T)
\text{ in }\Chowk[T,T^{-1}]$.
\end{rem}

\section{Products}
\label{abschnprodukte}
In this section, we investigate zeta functions of products of varieties whose
zeta functions satisfy a functional equation.

For the class of a Chow motive $x\in\Chowk$, we define
$Q^x(T):=\sum_{i\ge 0}\Alt^i(x)T^i$ and
$P^x(T):=\sum_{i\ge 0}\Sym^i(x)T^i$.

\begin{prop} \label{produkte}
The property of having a rational zeta function with values in $\Chowk$ 
satisfying a functional
equation is closed under products. More precisely, suppose that $X$ is an
$n$--dimensional smooth projective variety such that $[h(X)]=[X^+]+[X^-]$,
where $[\Alt^i(X^+)]=0$ for $i>e(X^+)$ and $[\Sym^i(X^-)]=0$ for $i>f(X^-)$.
Suppose furthermore that
$$Q^{X^+}(\frac{1}{\tate^n T})=T^{-e(X^+)}\tate^{-\frac{ne(X^+)}{2}}Q^{X^+}(T)
\text{ in }\Chowk[T,T^{-1}]$$
and 
$$P^{X^-}(\frac{1}{\tate^n
T})=T^{-f(X^-)}\tate^{-\frac{nf(X^-)}{2}}P^{X^-}(T)\text{ in }\Chowk[T,T^{-1}],$$
and likewise for a smooth projective variety $Y$. Then the same holds for
$X\times Y$, where $(X\times Y)^+=X^+\otimes Y^+\oplus X^-\otimes Y^-$,
$(X\times Y)^-=X^+\otimes Y^-\oplus X^-\otimes Y^+$, 
$e((X\times Y)^+)=e(X^+)e(Y^+)+f(X^-)f(Y^-)$ and 
$f((X\times Y)^-)=e(X^+)f(Y^-)+f(X^-)e(Y^+)$.
\end{prop}

We start with a special case.

\begin{lem}
Suppose that $x\in\Chowk$ is the class of a Chow motive satisfying
$$\deg Q^x=e \text{ and }Q^x(\frac{1}{\tate^m T})=
T^{-e}\tate^{-\frac{me}{2}}Q^x(T)$$
and that $y\in\Chowk$ is the class of a Chow motive satisfying
$$\deg Q^y=f\text{ and }
Q^y(\frac{1}{\tate^n T})=T^{-f}\tate^{-\frac{nf}{2}}Q^y(T).$$
Then the class $xy\in\Chowk$ satisfies
$$\deg Q^{xy}=ef\text{ and }Q^{xy}(\frac{1}{\tate^{m+n} T})=
T^{-ef}\tate^{-\frac{(m+n)ef}{2}}Q^{xy}(T).$$
\end{lem}

\begin{proof}[Proof of Lemma]
Denote the elementary symmetric functions in $\xi_1,\dots,\xi_e$ by 
$\sigma_i$ and the
elementary symmetric functions in $x_1,\dots,x_f$ by $s_i$.
Consider the following commutative diagram
\begin{equation*}\!\!\xymatrix@C=0.4cm{
\ganz[\xi_1,\xi_1^{-1}\!\!\!,\dots,\xi_e,\xi_e^{-1}\!\!\!,s_1,\dots,s_f,L,L^{-1}\!\!\!,t,t^{-1}]\ar[r]^-{\psi'}&
\Chowk[\xi_1,\xi_1^{-1}\!\!\!,\dots,\xi_e,\xi_e^{-1}\!\!\!,T,T^{-1}]\\
&\Chowk[\sigma_1,\dots,\sigma_e,T,T^{-1}]\ar@{^{(}->}[u]\ar[d]\\
\ganz[\sigma_1,\dots,\sigma_e,s_1,\dots,s_f,L,L^{-1}\!\!\!,t,t^{-1}]\ar^-{\varphi}[r]\ar[ru]
\ar[rd]\ar@{^{(}->}[uu]\ar@{_{(}->}[dd] &\Chowk[T,T^{-1}]\\ 
&\Chowk[s_1,\dots,s_f,T,T^{-1}]\ar[u]\ar@{_{(}->}[d]\\
\ganz[\sigma_1,\dots,\sigma_e,x_1,x_1^{-1}\!\!\!,\dots,x_f,x_f^{-1}\!\!\!,L,L^{-1}\!\!\!,t,t^{-1}]\ar^-{\psi''}[r]&
\Chowk[x_1,x_1^{-1}\!\!\!,\dots,x_f,x_f^{-1}\!\!\!,T,T^{-1}]}
\end{equation*}
where $\varphi(t)=T$, 
$\varphi(\sigma_i)=\Alt^i(x)$, $\varphi(s_j)=\Alt^j(y)$ and
$\varphi(L)=\tate$.

We know that 
\begin{align*}
q^x(t)&:=\prod_{1\le i\le e}(1+\xi_i t)\\
q^y(t)&:=\prod_{1\le j\le f}(1+x_j t)\\
q^{xy}(t)&:=\prod_{\substack{1\le i\le e\\1\le j\le f}}(1+\xi_i x_j t)\\
\end{align*}
are mapped by $\varphi$ to $Q^x(T)$, $Q^y(T)$ and $Q^{xy}(T)$,
and similarly, $q^x(\frac{1}{L^mt})$ is mapped to $Q^x(\frac{1}{\tate^mT})$,
and so on.

Now 
$$q^{xy}(\frac{1}{L^{m+n}t})=\prod_{1\le i\le e}q^y(\frac{\xi_i}{L^{m+n}t})
=\prod_{1\le i\le e}q^y(\frac{1}{L^nt_i}),$$
where $t_i=\frac{L^m t}{\xi_i}$.

We know that
$$\psi'(q^y(\frac{1}{L^nt_i}))=\psi'(t_i^{-f}L^{-\frac{nf}{2}}q^y(t_i))$$
and 
$$q^y(t_i)=\prod_{1\le j\le f}(1+x_j t_i)=t_i^f\prod_{1\le j\le f}x_j
\prod_{1\le j\le f}(1+\frac{1}{x_jt_i}).$$
Therefore,
$$\psi'(q^y(\frac{1}{L^nt_i}))=\psi'(L^{-\frac{nf}{2}}\prod_{1\le j\le f}x_j
\prod_{1\le j\le f}(1+\frac{\xi_i}{L^m x_jt}))$$
and hence
\begin{align*}
\varphi(q^{xy}(\frac{1}{L^{m+n}t}))
&=\varphi(L^{-\frac{nef}{2}}(\prod_{1\le j\le f}x_j)^e 
\prod_{1\le i\le e}\prod_{1\le j\le f}(1+\frac{\xi_i}{L^m x_jt}))\\
&= \varphi(L^{-\frac{nef}{2}}(\prod_{1\le j\le f}x_j)^e 
\prod_{1\le j\le f}q^x(\frac{1}{L^m \theta_j})),
\end{align*}
where $\theta_j=x_jt$.
As
$$\psi''(q^x(\frac{1}{L^m
  \theta_j}))=\psi''(\theta_j^{-e}L^{-\frac{me}{2}}q^x(\theta_j)),$$
we conclude that
\begin{align*}
\varphi(q^{xy}(\frac{1}{L^{m+n}t}))&=
\varphi(L^{-\frac{nef}{2}}
L^{-\frac{mef}{2}}t^{-ef}(\prod_{1\le j\le f}x_j)^e \prod_{1\le j\le f}x_j^{-e}q^x(\theta_j))\\
&=\varphi(t^{-ef}L^{-\frac{(m+n)ef}{2}}q^{xy}(t)).
\end{align*}
\end{proof}

Now due to Equations \ref{symsym}, \ref{symalt} and \ref{altsym} we have
the same behavior for $Q^{xy}(T)$ if 
$P^x(T)$ and $P^y(T)$ fulfill similar conditions,
and likewise for $P^{xy}(T)$ given the conditions for $P^x(T)$ and $Q^y(T)$ or
for $Q^x(T)$ and $P^y(T)$. This establishes the Proposition.

\end{document}